\documentclass[twocolumn]{autart}

\usepackage{graphicx}          
\usepackage[latin1]{inputenc}     

\usepackage{fancyhdr}                   
\usepackage{amsmath}                    
\usepackage{amsfonts}
\usepackage{amssymb}
\usepackage{mathrsfs}
\usepackage{stmaryrd}
\usepackage{subfigure}
\usepackage{pst-all,moredefs}
\usepackage[T1]{fontenc}   
\usepackage{dsfont}  
\usepackage{enumerate}
\usepackage{indentfirst}
\usepackage{lscape}
\usepackage{rotating}
\usepackage{epsfig,psfrag}
\usepackage{color,multicol}
\usepackage{calc,ifthen}

\newtheorem{rmrk}{Remark}[section]
\newtheorem{assm}{Assumption}[section]
\newtheorem{lemma}{Lemma}[section]
\newtheorem{proo}{Proof}[section]

\def\hal{{1 \over 2}}

\def\col{\mbox{col}}
\newcommand{\rank}{ \mbox{rank}}

\def\L2{{\cal L}_2}
\def\L2e{{\cal L}_{2e}}

\def\rea{\mathbb{R}}

\def\diag{\mbox{diag}}

\def\begequarr{\begin{eqnarray}}
\def\endequarr{\end{eqnarray}}
\def\begequarrs{\begin{eqnarray*}}
\def\endequarrs{\end{eqnarray*}}
\def\begarr{\begin{array}}
\def\endarr{\end{array}}
\def\begequ{\begin{equation}}
\def\endequ{\end{equation}}
\def\lab{\label}
\def\begdes{\begin{description}}
\def\enddes{\end{description}}
\def\begenu{\begin{enumerate}}
\def\begite{\begin{itemize}}
\def\endite{\end{itemize}}
\def\endenu{\end{enumerate}}

\def\lef[{\left[\begin{array}}
\def\rig]{\end{array}\right]}

\def\begcen{\begin{center}}
\def\endcen{\end{center}}
\def\texqq{\textbf{\texttt{q\hspace{-1.5mm}|}}}
\def\texqqs{\textbf{\texttt{q\hspace{-1mm}|}}}
\def\texpp{\textbf{\texttt{{|\hspace{-1.5mm}p}}}}
\def\texpps{\textbf{\texttt{{|\hspace{-1.2mm}p}}}}

\begin{document}

\begin{frontmatter}

\title{Two Globally Convergent Adaptive Speed  Observers for Mechanical Systems  } 
\thanks[footnoteinfo]{Corresponding author J.G. Romero Tel. +33 4 67 14 95 68 .}

\author[LIRMM]{Jose Guadalupe Romero* }\ead{Jose.Romero-Velazquez@lirmm.fr},    
\author[LSS]{Romeo Ortega}\ead{ortega@lss.supelec.fr},               

\address[LIRMM]{ Laboratoire d'Informatique, de Robotique et de Microélectronique de Montpellier, 161 Rue Ada, 34090 Montpellier, France}
\address[LSS]{ Laboratoire des Signaux et Syst\'emes, Sup\'elec, Plateau du Moulon, 91192 Gif-sur-Yvette, France}  
        
  
\begin{abstract}                          
A globally exponentially stable speed observer for mechanical systems was recently reported in the literature, under the assumptions of known (or no) Coulomb friction and no disturbances. In this note we
propose and adaptive version of this observer, which is robust {\em vis--\`a--vis} constant disturbances. Moreover, we propose a new globally convergent speed observer that, besides rejecting the disturbances,
estimates some unknown friction coefficients for a class of mechanical systems that contains several practical examples.
\end{abstract}

\end{frontmatter}

\section{Introduction}
\label{sec1}
The design of speed observers for mechanical systems is a problem of great practical importance that has attracted the attention of researchers for over $25$ years---the reader is referred to the recent books \cite{IIbook,besbook} for an exhaustive list of references. The first globally exponentially convergent speed observer for general, simple mechanical systems  was recently reported in \cite{ASTORTVEN}, where the Immersion and Invariance (I\&I) techniques developed in \cite{IIbook} were used. Although the observer of  \cite{ASTORTVEN} considers the case of systems with non--holonomic constraints, it relies on the assumptions of no friction and the absence of disturbances. In \cite{ROMORT} this speed observer was redesigned to accommodate the presence of {\em known} Coulomb friction and it was used to design a uniformly globally exponentially stable tracking controller using only position feedback for mechanical systems (without non--holonomic constraints).

In this paper we propose two new {\em robust} velocity observers for mechanical systems (without non--holonomic constraints). First, we add an adaptation stage to the observer of  \cite{ASTORTVEN} to reject constant input disturbances. Second, for mechanical systems with {\em zero Riemann symbols} (ZRS) \cite{NIJVAN,SPI,VENetal}, we propose a new adaptive speed observer that, besides rejecting the disturbances, estimates some of the {\em unknown friction coefficients}. It should be noted that in this case there are products of unmeasurable states and unknown parameters, a situation for which very few results are available in the observer design literature---even for the case of linear systems. {A similar transformation has been presented in  \cite{DO,DO1}, where a coordinate transformation that removes the quadratic terms in velocity is found to solve challenging position--feedback tracking problems for surface ships and mobile robots.}

The paper is organized as follows. The two adaptive observation problems addressed in the paper are presented in Section \ref{sec2}. The standing assumptions and a preliminary lemma are given in Section \ref{sec3}. In Section \ref{sec4} we present an adaptive observer for systems with unknown friction and disturbances. In Section \ref{sec5} the I\&I observer of \cite{ROMORT} is redesigned to accommodate the possible presence of disturbances and known friction. Some physical examples are given in Section \ref{sec6}. The paper is wrapped--up with some future research in Section \ref{sec7}. 

 \emph{Notation.} To avoid cluttering the notation, throughout the paper $\kappa$ and $\alpha$ are generic positive constants. $I_n$ is the $n \times n$ identity matrix and $0_{n \times s}$ is an
$n \times s$ matrix of zeros, $0_n$ is an $n$--dimensional column vector of zeros. Given $a_i \in \rea,\; i \in \bar n := \{1,\dots,n\}$, we denote with $\col(a_i)$ the $n$--dimensional column vector with
elements $a_i$. For any matrix $A \in \rea^{n \times n}$, $(A)_{i} \in \rea^n$ denotes the $i$--th column, $(A)^{i}$ the $i$--th row and $(A)_{ij}$ the $ij$--th element. That is, with $e_i\in \rea^n,\; i \in
\bar n$, the Euclidean basis vectors, $(A)_i := A e_i$, $(A)^{i}:=e_i^\top A$ and $(A)_{ij}:= e_i^{\top} A e_j$. For $x \in \rea^n$, $S \in \rea^{n \times n}$, $S=S^\top
>0$, we denote the Euclidean norm $|x|^2:=x^\top x$, and the weighted--norm $\|x\|^2_S:=x^\top S x$. Given a function $f:  \rea^n \to \rea$ we define the differential operators
$$
\nabla f:=\left(\frac{\displaystyle \partial f }{\displaystyle \partial x}\right)^\top,\;\nabla_{x_i} f:=\left(\frac{\displaystyle
\partial f }{\displaystyle \partial x_i}\right)^\top,
$$
where $x_i \in \rea^p$ is an element of the vector $x$. For a mapping $g : \rea^n \to \rea^m$, its Jacobian matrix is defined as
$$\nabla g:=\left [\begin{array}{cc}(\nabla g_1)^\top \\
\vdots\\ (\nabla g_m)^\top \end{array}\right],$$ where $g_i:\rea^n \to \rea$ is the $i$-th element of $g$.

\section{FORMULATION OF TWO ROBUST SPEED OBSERVATION PROBLEMS}
\label{sec2}
%
In the paper we consider $n$--degrees of freedom, {\em perturbed}, simple,  mechanical systems described in port--Hamiltonian (pH) form \cite{AvdS} by
\begin{equation}
\label{sys}
 \left[ \begin{array}{c} \dot{q}   \\
     \dot{\mathbf{p}} \end{array}\right]= \left[ \begin{array}{cc} 0 & I_n\\
     -I_n & -\mathbf{\mathfrak{R}}
    \end{array}\right]\nabla {H(q,\mathbf{p})} +\left[ \begin{array}{c} 0\   \\
     G(q) \end{array}\right]u+\left[ \begin{array}{c} 0\   \\
     d \end{array}\right]
\end{equation}
with total energy function $H:\rea^n \times \rea^n \to \rea$
\begin{equation}
\lab{h}
H(q,\mathbf{p})=\frac{1}{2}\mathbf{p}^\top M^{-1}(q)\mathbf{p}+V(q),
\end{equation}
where $q,\mathbf{p} \in \mathbb{R}^n$ are the generalized positions and momenta, respectively, $u \in \mathbb{R}^m$ is the control input, $G:\rea^n \to \rea^{n \times m}$ is the input matrix, the inertia
matrix $M:\rea^n \to \rea^{n \times n}$ verifies $M(q)=M^\top(q)
> 0$ and $V:\rea^n \to \rea$ is the potential energy function. As customary in the observer literature, it is assumed that
the control signal $u(t)$ is such that trajectories exist for all $t\geq 0$.

The system is subject to two different perturbations.
\begite
\item[-] \emph{Unknown} constant disturbances $d=\col(d_i) \in \rea^n$.
\item[-] Coulomb friction captured by
\begequ
\lab{frimat}
\mathbf{\mathfrak{R}}=\diag\{r_1,r_2,.,r_n\}\in \rea^{n \times n},
\endequ
with \emph{unknown} $r_i \geq 0,\; i \in \bar n$.
\endite

The problem is to design a globally convergent robust adaptive observer for the momenta $\mathbf{p}$. The main contributions of the paper are the following.
\begin{itemize}
\item[(i)] For systems with ZRS design a new observer that is globally convergent in spite of the presence of the disturbance $d$ and {\em some} unknown friction coefficients $r_i$.
\item[(ii)]  If the friction is known robustify the observer of  \cite{ASTORTVEN} to reject the disturbance $d$.\\
\end{itemize}
\begin{rmrk} \em
\lab{rem1}
The qualifier ``some" in item (i) is essential because, as will become clear later, except for the case of constant inertia matrix, we will not be able to consider the presence of unknown friction in all generalized
coordinates.
\end{rmrk}
\begin{rmrk} \em
\lab{rem2}
Notice that the objective is to observe only the momenta (equivalently, the velocity) not to ensure consistent estimation of the parameters $d$ and $r:=\col(r_i)$. As is well--known in the identification literature, a
necessary condition for parameter convergence is that the signals satisfy a persistency of excitation condition \cite{LJU}. In this respect, we notice that the system \eqref{sys}, \eqref{h}, which can be
represented in the state space form
$$
\begin{aligned}
\dot q & =  v \\
\dot v & =  U(q,v)  -\mathbf{\mathfrak{R}} v + G(q) u +d \\
\dot r & =  0\\
\dot d & = 0,
\end{aligned}
$$

for some $U:\rea^n \times \rea^n\to \rea^{n}$, with measurement $q$ does not satisfy the observability rank condition \cite{ISI} at zero velocity, hampering the observation of the parameters $r$.\footnote{The authors thank Prof.
Witold Respondek for this insightful remark.}
\end{rmrk}

\begin{rmrk} \em
See Remark 1 in \cite{ROMDONORT} for a physical interpretation of the disturbances $d$ that, we underscore, enter at the level of the momenta.
\end{rmrk}
%

\section{A Suitable pH Representation}
\label{sec3}
%
As shown in \cite{VENetal}, the change of coordinates
$$
(q,p) \mapsto(q, T^\top(q)\mathbf{p}),\hspace{3mm}
$$
where $T:\rea^n \rightarrow \rea^{n \times n}$ is a full rank factorization of the inverse inertia matrix, that is,
\begequ
\lab{mminonne}
M^{-1}(q)=T(q) T^\top(q),
\endequ
transforms \eqref{sys} into
\begin{align}
\nonumber
\left[ \begin{array}{c}
\dot{q} \\ \dot{p}
\end{array}\right] & = 
 \left[ \begin{array}{cc}
0 & T(q) \\ -T^\top(q) & J(q,p) -R(q)
\end{array}\right]
\nabla W(q,p) \\
& +  \left[ \begin{array}{cc} 0 \\  T^\top(q)G(q)
\end{array}\right] u + \left[ \begin{array}{cc} 0 \\ T^\top(q) d
\end{array}\right] ,
\label{sys_2}
\end{align}

with new Hamiltonian $W:\rea^n \times \rea^n \to \rea$
\begin{equation*}
\label{bar_H}
W(q,p)=\frac{1}{2}|p|^2+V(q),
\end{equation*}
the $jk$ element of the skew--symmetric matrix $J: \rea^n \times \rea^n \to \rea^{n \times n}$ given by
\begequ
\lab{jik}
{J}_{jk}(q,p) = -p^{\top}[(T)_{j}, (T)_{k}],
\endequ
with $[\cdot, \cdot]$ the standard Lie bracket \cite{SPI} and the transformed friction matrix
\begequ
\lab{r}
R(q):=T^\top(q)\mathbf{\mathfrak{R}}T(q) \geq 0.
\endequ

\begin{rmrk} \em
One possible choice of the factorization \eqref{mminonne} is the Cholesky factorization \cite{HORJOH,JAIROD}. But, as will become clear below, other choices may prove more suitable for the solution of the problem.
\end{rmrk}
%

\section{ROBUST OBSERVER FOR A CLASS SYSTEMS WITH FRICTION AND DISTURBANCES }
\label{sec4}
In this section we solve the problem of robust observation of momenta in the presence of disturbances and some unknown friction coefficients for a class of mechanical systems.
\subsection{Assumption on $M(q)$ and a preliminary lemma}
\label{subsec4.1}

\begin{assm} \em
\lab{ass1}
$M^{-1}(q)$ admits a factorization \eqref{mminonne} with a factor $T(q)$ verifying
\begequ
\lab{titj}
\Big[\Big(T(q)\Big)_{i}, \Big(T(q)\Big)_{j}\Big]=0,\;i,j \in \bar n.
\endequ
\end{assm}

Instrumental for the developments  of this paper is the following result, whose proof may be found in \cite{VENetal}.

\begin{lemma} \em
The following statements are equivalent:
\begenu
\item[(i)] $M(q)$ satisfies Assumption \ref{ass1}.
\item[(ii)] The Riemann symbols of $M(q)$ are all zero.\footnote{See equations (6) and (7) of \cite{VENetal} for the definition of these symbols.}
\item[(iii)] There exists a mapping $Q:\rea^n \to \rea^n$ such that
\begequ
\lab{qtminone}
\nabla Q(q) = T^{-1}(q).
\endequ
\endenu
\end{lemma}

\begin{rmrk} \em
Mechanical systems verifying Assumption \ref{ass1} have been extensively studied in analytical mechanics and have a deep geometric significance---stemming from Theorem 2.36 in \cite{NIJVAN}. They belong to the
class of systems that are partially linearizing via change of coordinates studied in \cite{VENetal}---see that paper for some additional references.
\end{rmrk}

\subsection{Assumptions on friction}
\label{subsec4.2}
To design our robust adaptive observer, besides Assumption \ref{ass1}, a restriction on the friction coefficients is imposed. Namely, we assume that there are $s$, with $s \leq n$, unknown coefficient and decompose the friction matrix $\mathbf{\mathfrak{R}}$ \eqref{frimat} as
$$
\mathbf{\mathfrak{R}}=\mathbf{\mathfrak{R}_k}+\mathbf{\mathfrak{R}_u}
$$
where $\mathbf{\mathfrak{R}_k},\mathbf{\mathfrak{R}_u}$ are $n \times n$ diagonal matrices containing the {\em known} and the {\em unknown} friction coefficients respectively. As a working example consider the case $n=3$ and $s=2$ with
$$
\mathbf{\mathfrak{R}_k} = \diag\{0,r_2,0\},\;\mathbf{\mathfrak{R}_u} = \diag\{r_1,0,r_3\}
$$
Similarly, with an obvious definition, we decompose the transformed friction matrix \eqref{r} into
$$
R(q)=R_k(q)+R_u(q).
$$
To streamline the presentation all friction coefficients are grouped in a vector $r=\col(r_i) \in \rea^n$ with  the unknown and known coefficients in  vectors $r_u \in \rea^s$ and $r_k \in \rea^{n-s}$, respectively. Thus, for our working example we have
$$
r=\col(r_1,r_2,r_3),\;r_k=r_2,\;r_u=\col(r_1,r_3).
$$
We also   define a set of integers $\kappa \subset \bar n$ that contains the indices of the unknown coefficients of $r$, which in the example is $\kappa=\{1,3\}$.

Finally, we define a matrix $C \in \rea^{n \times s}$  such that
\begin{equation}
\label{vecr}
C^\top r= r_u .
\end{equation}
Clearly, the matrix $C$ verifies:
\begite
\item $\rank\{C\}=s$.
\item For $j \in \kappa$, $(C)_j=e_{\kappa_j}$ .
\endite
In our example
$$
C = \lef[{ccc} 1 & 0 \\ 0 & 0 \\ 0 & 1 \rig].
$$
The following assumption and lemma are instrumental for our future developments.

\begin{assm}\em
\lab{ass2}
The $i$--th row of factor $T(q)$ is {\em independent of $q$} for $i \in \kappa$.\footnote{Consequently, the unknown friction coefficients are located in these rows.}   
\end{assm}

\begin{lemma} \em
Under Assumption \ref{ass2},  there exists  {\em constant} matrices $Y_j \in \rea^{n \times s},\;j \in \bar n,$ such that, for all vectors $z=\col(z_i)\in \rea^n$ we have
\begequ
\lab{linpar}
R_u(q)z=(\sum_{j=1}^n Y_j z_j)r_u.
\endequ
\end{lemma} 
\begin{proo} \em
From \eqref{vecr} it follows that 
\begin{equation}
\label{barrRu}
\mathfrak{R}_u=\sum_{i=1}^n e_i e_i^{\top}  (e_i^\top C r_u).
\end{equation}
Using the definition of $R_u(q)$ we get 
\begin{eqnarray}
 R_u(q)z& =&   T^\top (q)\left[  \sum_{i=1}^ne_i e_i^{\top}(e_i^\top C r_u)  \right] T(q)z,\nonumber  \\
 \nonumber
 &=&\sum_{i=1}^n T^\top (q)e_i e_i^{\top} T(q) z e_i^\top C r_u\\
 \label{Rux1} 
 &=&\sum_{i=1}^n \left[\sum_{j=1}^n T^\top (q)e_i e_i^{\top} T(q) \right]e_j z_j   e_i^\top Cr_u.
\end{eqnarray}
Hence, \eqref{linpar} follows swapping the sums and defining
\begin{eqnarray}
\lab{yj}
Y_j&:=&\sum_{i=1}^n T^\top (q)e_i e_i^{\top} T(q) e_j  e_i^\top C \nonumber   \\
&=&   \sum L_i e_j  e_i^\top C, \;j \in \bar n
\end{eqnarray}
 with  matrices $L_i$ defined as  
 \begequ
 \label{li}
L_i:= T^\top (q)e_i e_i^{\top} T(q) \hspace{1mm},\;i \in \bar n 
 \endequ
 
It only remains to prove that the matrices $Y_j$  and  $L_i$ are constant. Towards this end we refer to \eqref{yj} and notice that,  in view of Assumption \ref{ass2}, the term $e_i^\top T(q)$ is constant for $i \in \kappa$ while the term $e_i^\top C$ is an $1 \times s$ zero vector for $i \notin \kappa$. Completing the proof.
\end{proo}

\begin{rmrk} \em
As indicated in Remark \ref{rem1}, except for the case when  $M$ (and, consequently, the factor $T$) are constant, to satisfy Assumption \ref{ass2} we have to assume that some of the elements of $r$ are known. See Section \ref{sec6} for some physical examples.
\end{rmrk}

\subsection{First robust momenta observer}
\label{subsec4.3}
%
\begin{prop} \em
\label{pro1}
Consider  the system \eqref{sys}, \eqref{h} where the inertia matrix $M(q)$ and the friction matrix $\mathbf{\mathfrak{R}}$ verify Assumptions \ref{ass1} and \ref{ass2}. The $2n+s$ dimensional I$\&$I adaptive momenta observer
$$
\begin{aligned}
 \dot{p}_I & =   - T^\top(q)[\nabla V - G(q)u-\hat d]\\
 & -  (\sum_{i=1}^n Y_i \hat p_i)\hat r_u -[\lambda Q(q)+R_k(q)] \hat p\\
\dot r_{u_I} & =   (\sum_{i=1}^n Y^\top_i \hat p_i)(\dot p_I + \lambda \hat p)\\
\dot d_I & =  T(q) \hat p\\
\hat p & =  p_I+{\lambda}Q(q)\\
\hat r_u & =   r_{u_I} + {1 \over 2 \lambda}(\sum_{i=1}^s \hat p^\top L_i \hat p) e_i \\
\hat d & =   d_I+ q  \\
\hat {\mathbf{p}} & =  T^{-\top}(q)\hat p
\end{aligned}
$$
with the constant $n \times n$ matrices $L_i$ given by \eqref{li}, $Q(q)$ given in \eqref{qtminone}, $Y_i \in \rea^{n \times s},\;i \in \bar n$ given in \eqref{yj} and $\lambda>0$ a free parameter, ensures boundedness of all signals and
\begequ
\lab{momesterr}
\lim_{t\rightarrow \infty} [\hat {\mathbf{p}}(t)-\mathbf{p}(t)] =0.
\endequ
for all initial conditions $(q(0),\mathbf{p}(0)) \in \rea^n\times \rea^n.$
\end{prop}

%
\begin{proo} \em
Let the observation and parameter estimation errors be defined as
\begin{align}
\lab{z1}
\tilde p & = \hat p  - p \nonumber \\
 \tilde r_u& =  \hat r_u - r_u \\
 \tilde d &=  \hat d - d.
\end{align}
Following the I\&I adaptive observer procedure \cite{IIbook} we propose to generate the estimates as the sum of a proportional and an integral term, that is,
\begin{align}
\nonumber
\hat p&= p_I + p_P (q) \\
\hat r_u  &=  r_{u_I} + r_{u_P} (\hat p)\\
\hat d & =  d_I + d_P (q),
\lab{hatp0}
\end{align}
where the mappings
$$
\begin{aligned}
 p_P & :  \rea^n  \rightarrow \rea^n\\
 r_{u_P} & :  \rea^n\rightarrow \rea^{s} \\
 d_P & :  \rea^n\rightarrow \rea^{n},
\end{aligned}
$$
and the observer states  $p_I,d_I \in\rea^n$ and $r_{u_I} \in\rea^s$ will be defined below.\footnote{The reason for the particular selection of the arguments of the proportional terms will become clear below.}

First, we study the dynamic behavior of $\tilde p$ and compute
$$
\begin{aligned}
\dot{\tilde p}  &=  \dot{p}_I +  \nabla{p_P}T(q)p+ T^\top (q)[\nabla V -G(q)u] \\
        &  +[R_k(q)+R_u(q)](\hat p -\tilde p)-T^\top(q) (\hat d - \tilde d),
\end{aligned}
$$
where we have invoked Assumption \ref{ass1} that ensures---via \eqref{jik} and \eqref{titj}---that $J(q,p)=0$, and used \eqref{z1} to obtain the terms in the second row. Invoking \eqref{linpar} we can write
$$
R_u(q)\hat p=(\sum_{i=1}^n Y_i \hat p_i)r_u.
$$
Hence, proposing
\begin{align}
\nonumber
 \dot{p}_I & =  -\nabla{p_P}{ T}(q) \hat p - T^\top(q)[\nabla V - G(q)u]\\
 & -  \Big(\sum_{i=1}^n Y_i \hat p_i\Big)\hat r_u - R_k(q)\hat p+ T^\top(q) \hat d,
\label{dotpi}
\end{align}
yields
\begin{align}
\nonumber
\dot {\tilde p} & =  -[ R(q)  + \nabla_q p_P T(q)]\tilde p -(\sum_{i=1}^n Y_i \hat p_i)\tilde r_u + T^\top(q) \tilde d\\
 & =  -[ R(q)  + \lambda I_n]\tilde p  -(\sum_{i=1}^n Y_i \hat p_i)\tilde r_u + T^\top(q) \tilde d,
\label{tilpdyn}
\end{align}
where to obtain the second equations we have selected
$$
p_P(q)=\lambda Q(q),
$$
with $Q(q)$ given in \eqref{qtminone}.

Now, the time derivative of $ \tilde r_u$ is given as
$$
\begin{aligned}
\dot  {\tilde r}_u &= \dot {r}_{u_I} + \nabla r_{u_P} \dot {\hat p} \\
&= \dot {r}_{u_I} + \nabla r_{u_P}[ \dot p_I +\nabla p_P  T(q)p] \\
&= \dot {r}_{u_I}  + \nabla r_{u_P}[ \dot p_I + \lambda(\hat p -\tilde p)].
\end{aligned}
$$
Hence, choosing
\begin{equation}\nonumber
\dot {r}_{u_I}= - \nabla r_{u_P} (  \dot p_I +\lambda \hat p ),
\end{equation}
yields
\begin{equation}
 \dot{\tilde r}_u= -\lambda \nabla r_{u_P}\tilde p.
\label{dotzax}
\end{equation}
Finally, the time derivative of $ \tilde d$ is given as
$$
\begin{aligned}
\dot  {\tilde d} &= \dot {d}_I + \nabla d_P T(q) p \\
&= \dot d_I + \nabla d_P T(q)(\hat p - \tilde p).
\end{aligned}
$$
Hence, choosing
\begin{equation}
\nonumber
\dot d_I= - \nabla d_P T(q) \hat p,
\end{equation}
yields
\begin{equation}
 \dot{\tilde d}= -\nabla d_P T(q)\tilde p.
\label{dottild}
\end{equation}
We will now analyze the stability of the  {\em error model} \eqref{tilpdyn}, \eqref{dotzax} and \eqref{dottild} with the aid of the proper Lyapunov function candidate
\begequ
\label{V21}
V(\tilde p, \tilde d, \tilde r_u)= {1\over 2}( |\tilde p|^2 +|\tilde d|^2 +|\tilde r_u|^2).
\endequ
Taking its time-derivative  we obtain
\begin{align}
\nonumber
\dot V&= -\tilde p^\top \Big[ R (q)+ \lambda I_n \Big]\tilde p -\tilde p^\top\big[(\sum_{i=1}^n Y_i \hat p_i)\tilde r_u - T^\top(q) \tilde d \big]\\
&  - \big[\lambda\tilde r_u^\top \nabla r_{u_P} +\tilde d^\top \nabla d_P T(q)\big] \tilde p.
\lab{dotv}
\end{align}
Clearly, {\em if} the mappings $r_{u_P}(\hat p)$ and  $d_P(q)$ solve the partial differential equations (PDEs)
\begin{equation}
\begin{aligned}
 \nabla r_{u_P} & =  -{1 \over \lambda} (\sum_{i=1}^n Y^\top_i \hat p_i)\\
 \nabla d_{P} & =  I_n,
 \end{aligned}
\label{PDE_betaa}
\end{equation}
one gets
\begequ
\label{dotV}
\dot V =-\tilde p^\top [ R (q)+ \lambda I_n] \tilde p \leq - \lambda |\tilde p|^2.
\endequ
From \eqref{V21}, \eqref{dotV} we conclude that $\tilde p$ $\in$ $\mathcal{L}_2$ $\cap$ $\mathcal{L}_{\infty}$ and $\tilde d, \tilde r_u \in \mathcal{L}_{\infty}$. Doing some standard signal chasing it is
straightforward to prove from here that \eqref{momesterr} holds.

Motivated by the conclusion above let us now study the PDEs \eqref{PDE_betaa}. The second one has the trivial solution $d_P=q$. Regarding the first one, it is clear that the elements of the mapping
$r_{u_P}(\hat p)$ must be of the quadratic form
$$
(r_{u_P}(\hat p))_i={1\over  2{\lambda}}  \hat p^ \top L_i \hat p,\;i \in \bar s,
$$
with {\em constant, symmetric} matrices  $L_i \in \rea^{n \times n}$. Replacing the expression above in the PDE \eqref{PDE_betaa} yields
$$
\begin{aligned}
\left[ \begin{array}{c}
 \hat p^\top L_1 \\
 \vdots\\
 \hat p^\top  L_{s}
\end{array}\right]= -\sum_{i=1}^n Y^\top_i\hat p_i.
\end{aligned}
$$
That lead us to the solution
\begequ
\lab{solpde}
 \left[ \begin{array}{ccc}
  L_1^\top e_j &  \hdots & L_{s}^\top e_j \end{array}\right]=-Y_j,\;j \in \bar n.
\endequ

It only remains to show that the resulting matrices $L_i$ are symmetric. From \eqref{solpde} we get
\begequarrs
-L_j^\top & = & -\lef[{cccc} L_j^\top e_1 &  L_j^\top e_2 & \hdots & L_j^\top e_n\rig]\\
                & = &  \lef[{cccc} Y_1 e_j & Y_2 e_j & \hdots & Y_n e_j\rig]
\endequarrs
Clearly, the matrix $L_j$ is symmetric if and only if
$$
e_i^\top Y_k = e_k^\top Y_i,\;\forall i,k \in \bar n,\; i \neq k.
$$
This fact can be easily verified using \eqref{yj} 
\begequarrs
e_k^\top Y_i & = & e_k^\top \sum_{j=1}^n T^\top (q)e_j e_j^{\top} T(q) e_i  e_j^\top C \\
 & = & e_i^\top \sum_{j=1}^n T^\top (q)e_j e_j^{\top} T(q) e_k  e_j^\top C \\
 & = &e_i^\top Y_k.
\endequarrs

Replacing all the derivations above in $\dot p_I$, $\dot d_I$ and $\dot r_{u_I}$ gives the equations given in the proposition completing the proof.
\end{proo}
%
\begin{rmrk} \em
If Assumption \ref{ass1} is not imposed a term $J(q,p)p$ appears in the error equation \eqref{tilpdyn} and \eqref{dotzax}. Even though this term is quadratic in the unknown state $p$, the properties of
$J(q,p)$ can be used to handle this term in the first error equations---this is done in the second observer in the next section. However, there is no obvious way to create a suitable error term for the second
error equation styming the relaxation of Assumption \ref{ass1}.
\end{rmrk}

\begin{rmrk} \em
If the matrices $Y_i$ are not constant the first PDE in \eqref{PDE_betaa} does not admit a solution, hence Assumption \ref{ass2} is required. It can be shown that making $r_{u_P}$ function of $q$ does not
solve the problem, because a quadratic function of $\hat p$ will appear in $\dot V$.
\end{rmrk}

\begin{rmrk} \em
As shown in Proposition 6 of \cite{VENetal} the dynamics of mechanical systems satisfying Assumption \ref{ass1} expressed in the coordinates $(Q,p)$ take the form
\begin{align}
\nonumber
\dot Q & =  p \\
\lab{douintsys}
\dot p & =  -\tilde R(Q)p-\tilde T^{\top}(Q)[\nabla \tilde V(Q) - \tilde G(Q) u+  d],
\end{align}
where $\tilde {(\cdot)}(Q):=(\cdot)(Q^I(Q))$, with $Q^I:\rea^n \to \rea^n$ a left inverse of $Q(q)$, that is, $Q(Q^I(z))=z$ for all $z \in \rea^n$. Although the construction of an observer for
\eqref{douintsys} when the friction is known is straightforward, the case of unknown friction is far from trivial. Applying the I$\&$I procedure used in Proposition \ref{pro1} leads, for the definition of
$d_P(q)$, to a PDE of the form
$$
\nabla S(Q) =  \tilde T(Q),
$$
whose solution is not obvious.
\end{rmrk}

\section{ROBUST OBSERVER FOR GENERAL PERTURBED MECHANICAL SYSTEMS WITH KNOWN FRICTION}
\lab{sec5}
%
In this section we redesign the I\&I speed observer of  \cite{ROMORT}, see also \cite{ASTORTVEN}, to ensure its global convergence in spite of the presence of the {\em unknown} disturbances $d$ and {\em known} friction forces in all coordinates.
\begin{prop} \em
\label{pro3}
Consider  the system  \eqref{sys}, \eqref{h} with {\em known} friction matrix $\mathfrak{R}$. There exist smooth mappings
$$
\begin{aligned}
\mathbf{A}&:\rea^{4n} \times \rea_{\geq 0} \times \rea^n\times \rea^n \to \rea^{4n+1}\\
\mathbf{B}&: \rea^{4n} \times \rea_{\geq 0} \times \rea^n\rightarrow \rea^{n}
\end{aligned}
$$
such that the interconnection of  \eqref{sys}, \eqref{h} with
\begin{equation}
\begin{aligned}
\dot{\mathrm{X}} &=  \mathbf{A}(\mathrm{X},q,u)\\
\label{state_eta}
\hat {\mathbf{p}}&=\mathbf{B}(\mathrm{X},q) ,
\end{aligned}
\end{equation}
where $\mathrm{X} \in  \rea^{4n} \times \rea_{\geq 0},\;\hat {\mathbf{p}} \in \rea^n$, ensures \eqref{momesterr} holds for all initial conditions
$$
(q(0), \mathbf{p}(0), \mathrm{X} (0)) \in \rea^n\times \rea^n  \times \rea^{4n} \times \rea_{\geq 0}.
$$
This implies that, in spite of the presence of the {\em unknown disturbances} $d$, \eqref{state_eta} is a globally convergent momenta observer for the mechanical system with friction  \eqref{sys}, \eqref{h}.\\
\end{prop}

\begin{proo} \em
The construction of the observer follows very closely the one reported in \cite{ROMORT} with the only difference of the inclusion of an adaptation law for the unknown disturbance parameters $d$. However, for the sake of completeness, a detailed derivation of all the steps is given.

Define the estimation errors
\begin{equation}
\begin{aligned}
\lab{z}
\tilde p & = \hat p  - p \\
 \tilde d &= \hat d- d.
 \end{aligned}
\end{equation}
Following the I\&I adaptive observer procedure \cite{IIbook} we propose to generate the estimates as
\begin{equation}
\begin{aligned}
\lab{hatp}
\hat p&:=& p_I + p_P (q,  \texqq, \texpp)  \\
\hat d&:=& {d_I + d_P (q,r)}
 \end{aligned}
\end{equation}

where the mappings $p_P:\rea^n \times \rea^n \times \rea^n \rightarrow \rea^n$ and $d_P \in \rea^n $, and the observer states $d_I \in \rea^{n}$,  $p_I \in\rea^n$ and $r \in \rea$ are defined such that
\eqref{momesterr} holds.

We, therefore, study the dynamic behavior of $\tilde p$ and compute
$$
\begin{aligned}
\dot {\tilde p}  &=  \dot{p_I} +  \nabla_{q}{p_P}{ \dot{q}} + \nabla_{\texqqs}{ p_P} \dot{\texqq} + \nabla_{\texpps}{ p_P}{ \dot \texpp}-\\
        &- J(q,p)p +  T^\top (q)[\nabla V -G(q)u] +R(q)p-T^\top(q)d.
\end{aligned}
$$
In \cite{ASTORTVEN} it has been shown that the mapping $J(q,p)$ defined in \eqref{jik} verifies the following properties:
\begite
\item[(P.i)] $J(q,p)$ is linear in the second argument, that is
 $$
 J(q, \alpha_1 p+\alpha_2 \bar{p}) = \alpha_1{J}(q,p) + \alpha_2{J}(q, \bar{p})
 $$

  for all $q$, $p$, $\bar{p}$ $\in \rea^n$, and  $\alpha_1$, $\alpha_2$ $\in \rea$.
\item[(P.ii)] There exists a mapping  ${\bar J} :\rea^{n} \times\rea^{n}\to \rea^{n\times n}$ satisfying
$${J}(q, p)\bar{p}= {\bar J}(q, \bar{p})p.$$
\endite
Hence, proposing
\begin{align}
\nonumber \dot{p_I} &:=  - \nabla_{\texqqs}{p_P}\dot{\texqq}- \nabla_{ \texpps}{p_P}{ \dot{\texpp}} + J(q,\hat p)\hat p -R(q)\hat p - \\
&- T(q)^\top(q)\nabla V + v-\nabla_{q}{p_P}T(q) \hat p  +T^\top(q) \hat d, \nonumber \\
\label{dot_xi}
\end{align}
together with Properties (P.i) and (P.ii) yields
\begequ
\label{zdyn}
\dot {\tilde p} = [J(q,p)  + {\bar J}(q, \hat p) -R(q)- \nabla_q p_P T(q)] \tilde p +T^\top(q)\tilde d.
\endequ
It is clear that if the mapping $p_P$ solves the PDE
$$
\nabla_q p_P = [\psi{ I_n} + {\bar J}(q, \hat p) ]T^{-1}(q),
$$
with $\psi >0$ a design constant, the $\tilde p$--dynamics reduces to
$$
\dot {\tilde p }= [J(q,p) - \psi { I_n}-R(q)] \tilde p  +T^\top(q) \tilde d.
$$
Recalling that $J(q,p)$ is skew--symmetric and $R(q) \geq 0$ the unperturbed part of the error dynamics above, {\em i.e.} when $\tilde d=0$, is exponentially stable.

Similarly to \cite{ROMORT}, to avoid the solution of the PDE, the dynamic scaling technique is used. Towards this end, define the mapping
\begequ
\label{ideal}
{\mathcal H}(q, \hat p):= [\psi { I_n} + {\bar J}(q, \hat p) ]T^{-1}(q).
\endequ
and define $p_P$ as
\begequ
\label{beta}
p_P (q, \texqq, \texpp):= {\mathcal H}({\texqq, \texpp})q.
\endequ
The choice above yields $\nabla_q p_P =  {\mathcal H}({\texqq, \texpp})$, which may be written as
\begequ
\nabla_q p_P = {\mathcal H}(q, \hat p) - [{\mathcal H}(q, \hat p) -{\mathcal H}({\texqq, \texpp})].
\label{nablabeta}
\endequ
Now, since the term in brackets in  (\ref{nablabeta}) is equal to zero if $\texpp=\hat p$ and $\texqq=q$, there exist mappings
$$
{\Delta}_q,{\Delta}_p :\rea^{n}\times\rea^{n}\times\rea^{n}\to\rea^{n\times n}
$$
verifying
\begequ
\label{Deltazero}
{ \Delta}_q(q,\texpp, 0)=0, \quad { \Delta}_p(q,\texpp, 0)=0,
\endequ
and such that
\begequ
\lab{deldel}
{\mathcal H}(q, \hat p) -{\mathcal H}({\texqq, \texpp})={ \Delta}_q(q, \texqq, e_q)  + { \Delta}_p(q, \texpp,  e_p),
\endequ
where
\begequ
\label{errors}
e_q := \texqq - q, \quad e_p := \texpp - \hat p.
\endequ

Substituting (\ref{ideal}), (\ref{nablabeta}) and \eqref{deldel} in (\ref{zdyn}), yields
$$
\begin{aligned}
\dot {\tilde p}& = [J(q,p)- \psi { I_n}-R] \tilde p \\
&  + \Big({\Delta}_q(q, \texqq, e_q)  +  {\Delta}_p(q, \texpp,  e_p)\Big)T(q) \tilde p +T^\top(q) \tilde d.
\end{aligned}
$$
The mappings ${\Delta}_q$, ${\Delta}_p$ play the role of disturbances that are dominated with a dynamic scaling and a proper choice of the observer dynamics. For, define the dynamically scaled
off--the--manifold coordinate
\begequ
\lab{eta}
\eta=\frac{\displaystyle 1}{\displaystyle r} \tilde p,
\endequ
where $r$ is a  scaling factor to be defined. The dynamic behavior of $\eta$ is given by
\begequ
\label{deta}
\dot{ \eta} = ({ J}-R-\psi{ I}){ \eta} + ({ \Delta_q}+{\Delta_p})T {\eta}  +\frac{1}{r}T^\top \tilde d- \frac{\dot r}{r}{ \eta},
\endequ
where, for brevity, the arguments of the mappings are omitted.

Mimicking \cite{ASTORTVEN} select the dynamics of $\texqq$, $\texpp$ as
\begin{equation}
\begin{aligned}
\label{dot_qq}
 { \dot \texqq} &= T(q)\hat p - \psi_1{ e_q}   \\
 \dot \texpp &= - T^\top(q)\nabla V + v + J(q, \hat p)\hat p -R\hat p\\
 & - \psi_2{ e_p}  +T^\top(q) \hat d
\end{aligned}
\end{equation}
where $\psi_1,\psi_2$ are some positive functions of the state defined later. Using \eqref{dot_qq}, together with (\ref{errors}), we get
\begin{equation}
\begin{aligned}
\label{doterrors}
 \dot{e}_q &= T (q) \eta r - \psi_1{e_q}   \\
 \dot{e}_p &= \nabla_{q}p_P T(q) \eta r - \psi_2{ e_p}.
\end{aligned}
\end{equation}
Moreover, select the dynamics of $r$ as
\begequ
\lab{dotr}
\dot r  =  -{\psi \over4}(r-1)+{r\over \psi }(\|{ \Delta_p T}\|^2+\|{ \Delta_q T}\|^2), \; r(0)\ge 1,
\endequ
with $\|\cdot\|$ the matrix induced $2$--norm. At this point we make the important observation that the set
$$
\{r\in\rea: r\geq 1\}
$$
is {\em invariant} for the dynamics (\ref{dotr}). Hence, $r(t) \geq 1,\;\forall t \geq 0$.

On other hand, taking the time-derivative of $\tilde d$, we get
$$
\begin{aligned}
\dot{\tilde d}&=\dot{d}_I+ \nabla_q d_PT(q) p + \nabla_r d_P\dot r \nonumber \\
&=\dot{d}_I+\nabla_q d_P(q,r) T(q)(\hat p-r\eta ) + \nabla_r d_P\dot r
\end{aligned}
$$
and choosing
\begequ
\lab{dotdi}
\dot d_I= -\nabla_q d_PT(q) \hat p - \nabla_r d_P \dot r,
\endequ
the $\tilde d$--dynamics take the form
 \begin{equation}
\dot{\tilde d}= - r\nabla_q d_P T(q) \eta
\label{dotza1}
 \end{equation}
 
We now analyze the  {\em error} system \eqref{eta}, \eqref{doterrors}, \eqref{dotr}, \eqref{dotza1}---with the  coordinate $\tilde r = (r-1)$. For, define the proper Lyapunov function candidate.
\begequ
\label{V2}
V(\eta,e_q,e_p,\tilde r, z_a):= \hal\left(|\eta|^2+|e_q|^2 + |e_p|^2+\tilde r ^2 +|\tilde d|^2 \right).
\endequ
Taking its time-derivative  we obtain
$$
\begin{aligned}
\dot V&\leq -\left({\psi \over 4} -1\right)|\eta|^2 - \left(\psi_1 - \hal r^2\|T(q)\|^2\right)|e_q|^2 \\
&   -\left(\psi_2 - \hal r^2\|\nabla_q p_P\|^2\|T(q)\|^2\right)| e_p|^2 + \tilde r \dot r+ \\
& + \tilde d^\top\left(\frac{1}{r}- r\nabla_q d_P\right)T(q) \eta 
\end{aligned}
$$
Clearly, if we set
\begequ
\lab{psipsi}
{\psi=4(1 + \psi_3)},\;\psi_1= \hal r^2\|T(q) \|^2 + \psi_4
\endequ
and
$$
\psi_2= \hal r^2\|\nabla_q p_P\|^2\|T(q) \|^2 + \psi_5,
$$
where  $\psi_3, \psi_4, \psi_5$ are positive functions of the state defined below, one gets
$$
\begin{aligned}
\dot V&\leq -\psi_3  |\eta|^2 - \psi_4|e_q|^2 - \psi_5|e_p|^2 + \tilde r \dot r \\
&+ \tilde d^\top\left(\frac{1}{r}- r\nabla_q d_P(q,r)\right)T(q) \eta .
\end{aligned}
$$
To eliminate the cross term appearing in the last term of the right hand side above we select
\begin{equation}
\label{d_P}
d_P(q,r)={1 \over r^2} q,
\end{equation}
which clearly solves the PDE
$$
{1 \over r} -r \nabla_q d_P(q,r)=0.
$$
It only remains to study the term $\tilde r \dot r$, which is given by
$$
\tilde r \dot r= -{\psi \over4}\tilde r^2+\tilde r {r\over \psi }(\|{ \Delta_p T}\|^2+\|{ \Delta_q T}\|^2).
$$
Now, (\ref{Deltazero}) ensures the existence of mappings $\bar \Delta_q$, $\bar \Delta_p:\rea^{n}\times\rea^{n}\times\rea^{n}\to\rea^{n\times n}$ such that
$$
\begin{aligned}
\|\Delta_q(q, \texpp, e_q)\| &\leq \|\bar \Delta_q(q, \texpp, e_q)\|\;|e_q| \\
\|\Delta_p(q, \texpp, e_p)\| &\leq \|\bar\Delta_p(q, \texpp, e_p)\|\;|e_p|.
\end{aligned}
$$
Hence
$$
\|{\Delta_p T}\|^2+\|{ \Delta_q T}\|^2\leq \|T\|^2|(\|\bar\Delta_p\|^2|e_p|^2 + |\bar\Delta_q\|^2|e_q|^2).
$$
Setting
$$
\begin{aligned}
\psi_3&= \kappa \\
 \psi_4&= {r\tilde r \over 4(1 + \psi_3)}\|T \|^2\|\bar\Delta_q\|^2 + \kappa\\
\psi_5&= {r\tilde r \over 4(1 + \psi_3)}\|T \|^2\|\bar\Delta_p\|^2 + \kappa,
\end{aligned}
$$
one  gets
\begin{equation}
\label{dotVt}
\dot V\leq - \kappa (|\eta|^2 + |e_q|^2 +|e_p|^2 +\tilde r^2]
\end{equation}
for some positive constant $\kappa$.

The proof is completed invoking the arguments of \cite{ROMORT}, selecting the observer state as
$$
{\mathrm{X}}:= (\texqq,\texpp,p_I,d_I,r-1),
$$
and defining $\mathbf{A}(\mathrm{X},q,u)$ from \eqref{dot_xi}, \eqref{dot_qq}, \eqref{dotr} and \eqref{dotdi} and setting $\mathbf{B}(\mathrm{X},q)$ via
\eqref{hatp}.
\end{proo}

\begin{rmrk} \em
It is clear from the proof that the key step to reject the disturbances is to make the proportional term of the parameter estimator, $d_P$ a function of the dynamic scaling factor $r$, see (\ref{d_P}).
\end{rmrk}

\section{PHYSICAL EXAMPLES}
\label{sec6}
%
In this section we present  three physical mechanical systems that satisfy the conditions of Proposition \ref{pro1}. Consequently, robust adaptive speed observation is possible for them.
%
\subsection{Constant inertia matrix}
In the case of constant inertia matrix Assumption \ref{ass1} is trivially satisfied, because the factor $T$ can be taken to be constant. Assumption \ref{ass2} is also satisfied with $\dim(q)=n$ and $C=I_n$, hence all friction coefficients can be identified.

Given any constant factor $T$, the vector field $Q(q)$ that solves  \eqref{qtminone} is given by
$$
Q(q)= T^{-1} q.
$$
Finally, from \eqref{yj} and \eqref{li} we get
$$
\begin{aligned}
Y_j&=\sum_{i=1}^n T^\top (q)e_i e_i^{\top} T(q) e_j  e_i^\top\;j \in \bar n \\
L_i &= [(T)^i]^\top (T)^i,\;i \in \bar n.
\end{aligned}
$$
\subsection{Planar redundant manipulator with one elastic degree of freedom}
\begin{figure}[h]
 \centering
\includegraphics[width=0.65\linewidth]{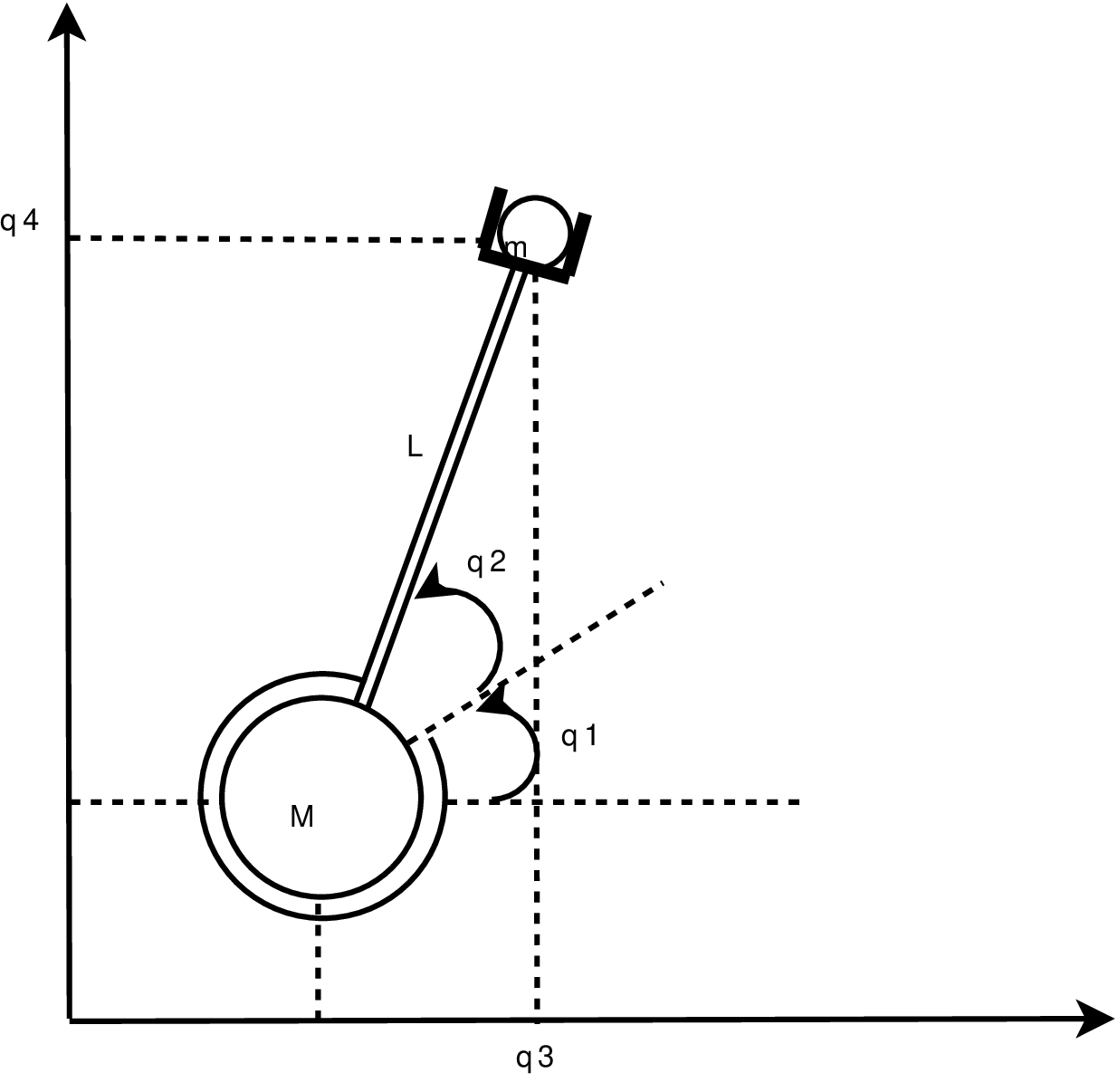}
\caption{Planar redundant manipulator with one elastic degree of freedom.}
   \label{fig7}
\end{figure}

This is a 4-dof underactuated mechanical system depicted in \ref{fig7}.
%
The inverse inertia matrix is given by
$$
M^{-1}(q)=\left[ \begin{array}{cccc} 1\over I & - 1 \over I &  0 & 0 \\
     * &  {1\over a_2^2}+ {1 \over I} & - {1 \over m\ell }S_{12} & {1 \over m\ell } C_{12}\\
     *& *& {1\over a_3^2}& 0 \\
      *& *& * & {1\over a_3^2}
    \end{array}\right],
$$
where we defined
$$
\begin{aligned}
S_{12}&:= \sin(q_1+q_2),\;C_{12}:=\cos(q_1+q_2)\\
a_2&:= {{\sqrt{Mm}\ell } \over \sqrt{m+M}},\;a_3:= \sqrt{M+m},
\end{aligned}
$$
and the definition of all constants may be found in \cite{VENetal}. The Cholesky factorization is given as
$$
T(q)=\left[ \begin{array}{cccc} 1\over   \sqrt I & 0 &  0 & 0 \\
     - 1 \over \sqrt I & {1 \over a_2}  & 0 & 0 \\
     0& - \sqrt{ Mï¿½\over m} {1\over a_3}S_{12}& 1\over a_3 & 0 \\
      0&  \sqrt{ Mï¿½\over m} {1\over a_3}C_{12}& 0 & 1\over a_3
    \end{array}\right],
$$
and the vector field $Q(q)$ that solves  \eqref{qtminone} is
$$
Q(q)= \left[ \begin{array}{c}
 \sqrt I q_1 \\
 a_2(q_1 +q_2)  \\
-  a_2 \sqrt{ Mï¿½\over m}  C_{12} +a_3 q_3 \\
-  a_2 \sqrt{ Mï¿½\over m}S_{12}+ a_3 q_4
    \end{array}\right].
$$
A matrix $C$ that satisfies Assumption \ref{ass2} is
$$
C = \lef[{cc} I_2 \\ 0_{2 \times 2} \rig].
$$
Hence, we can consider  as unknown the frictions in the elastic coordinate $r_1$ and the revolute joint $r_2$.

Finally,
$$
\begin{aligned}
Y_1^\top &= \left[ \begin{array}{cccc}
 1 \over I& 0& 0& 0   \\
 {1 \over I} &    -{a_2\over \sqrt I} & 0 &0     \end{array}\right]\\
Y_2^\top & =  \left[ \begin{array}{cccc}
0 & 0& 0& 0   \\
 -{a_2\over \sqrt I}  &   a_2^2 & 0 &0     \end{array}\right]\\
 L_1 & =   \left[ \begin{array}{ccc}
 \left[ \begin{array}{cc}
 {1 \over I} & 0\\
 0 &0
 \end{array}\right] & 0_{2 \times 2}   \\
 0_{2 \times 2}& 0_{2 \times 2}      \end{array}\right] \\
 L_2 & =   \left[ \begin{array}{ccc}
 \left[ \begin{array}{cc}
 {1 \over I} &  -{a_2\over \sqrt I} \\
 -{a_2\over \sqrt I}  & a_2^2
 \end{array}\right] & 0_{2 \times 2}   \\
 0_{2 \times 2}& 0_{2 \times 2}      \end{array}\right]
 \end{aligned}
 $$
\subsection{2D-Spider crane gantry cart}
\begin{figure}[h]
 \centering
\includegraphics[width=0.75\linewidth]{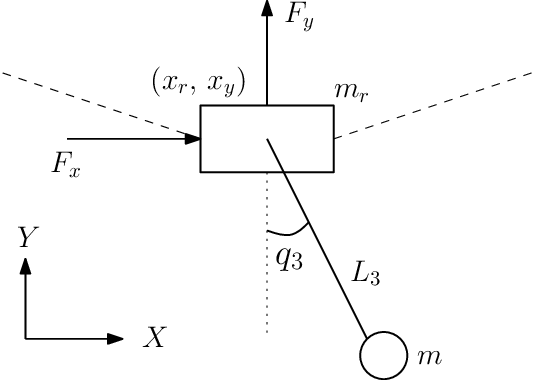}
 \caption{2D-Spider crane gantry cart}
 \label{fig1}
\end{figure}

\lab{subsec6.2}
This is a 3-dof underactuated mechanical system depicted in Fig. \ref{fig1}.
The inertia matrix is
$$
M(q)= \left[ \begin{array}{ccc} m_r+m & 0 &  {m L_3 }C_3 \\
     *& m_r+m & mL_3S_3  \\
     *& *&m L_3^2
         \end{array}\right],
$$
with inverse
$$
M^{-1}(q)= \left[ \begin{array}{ccc} {m_r+mC_3^2}\over{(m_r+m)m_r} & {mC_3S_3}\over{(m_r+m)m_r} &  -{C_3} \over{L_3 m_r} \\
     *& {{m_r+m-m C_3^2}\over {(m_r+m)m_r}} & -{{S_3}\over{m_r L_3} }  \\
     *& *&{ {m_r+m}\over {m_r L_3^2m}}
    \end{array}\right]
$$
where, to simplify the notation, we have defined
$$
S_3 := \sin(q_3),C_3:=\cos(q_3),
$$
and the definition of all constants may be found in \cite{KABAMULBO}. An upper triangular factorization of $ M^{-1}(q)$ is given as
\begequ
\lab{tt}
T(q)= \left[ \begin{array}{ccc} a & 0  &  -b C_3 \\
     0&  a & -b S_3 \\
     0& 0 & c
    \end{array}\right],
    \endequ
where we defined the constants
$$
a := {{1}\over{\sqrt{m_r+m}}}, b:={{1}\over{c L_3 m_r }},c:=\sqrt{{{m_r+m}\over{mL_3^2m_r}}}.
$$
We can  check that the columns of $T(q)$ satisfy \eqref{titj} and thus the system verifies Assumption 1. Taking the inverse of  $T(q)$ we get
$$
T^{-1}(q)= \left[ \begin{array}{ccc} {1}\over{a} & 0  &  aL_3mC_3 \\
   0&   {1}\over{a} & aL_3mS_3 \\
     0& 0 &  {1}\over{c}
    \end{array}\right],
$$
From the equation above it is clear that a mapping $Q(q)$ that solves  \eqref{qtminone}  is
$$
Q(q)=\left[ \begin{array}{c}
{1\over a} q_1 +a L_3 m S_3\\
{1\over a} q_2-a L_3 m C_3  \\
{1 \over c }q_3  \\
    \end{array}\right].
$$

From the definition of $T(q)$ in \eqref{tt} it is clear that
$$
C^\top = \lef[{ccc} 0 & 0 & 1\rig]
$$
satisfies Assumption \ref{ass2}. Hence, we can consider  as unknown the friction parameter $r_3$.

Finally, from \eqref{yj} and \eqref{li} we get
$$
\begin{aligned}
Y_1^\top  &= \left[ \begin{array}{ccc}
0 &0 & c^2 \\
    \end{array}\right],\\
L_1&=\left[ \begin{array}{ccc}
0 &0 &0 \\
0 &0 & 0 \\
0 &0 & c^2 \\
    \end{array}\right].
 \end{aligned}
 $$
We show several simulations of the proposed observer.  The system has the two forces shown in  Fig. \ref{fig1} as control inputs, that is,  $u=\col(F_x,F_y)$ and constant input matrix $G$ of the form
$$
G = \lef[{cc} 1 & 0 \\
                             0& 1\\
                             0 & 0\rig].
$$
The parameters are taken as $m_r=0.5$ $kg$ for ring mass, $m= 1$ $kg$ for payload mass and $L_3=0.5$ $m$ for the cable length.  We fix the control inputs as $F_x=1.535 \cos(t)$ and $F_y=7.67 \sin(t)$.  The disturbances are taken as $d=\col( 0.1,0.2,0.2)$ and the friction coefficients $r=\col(0,0,0.5)$, with $r_3$ being unknown. 

\begin{figure}[h]
 \centering
\includegraphics[width=1\linewidth]{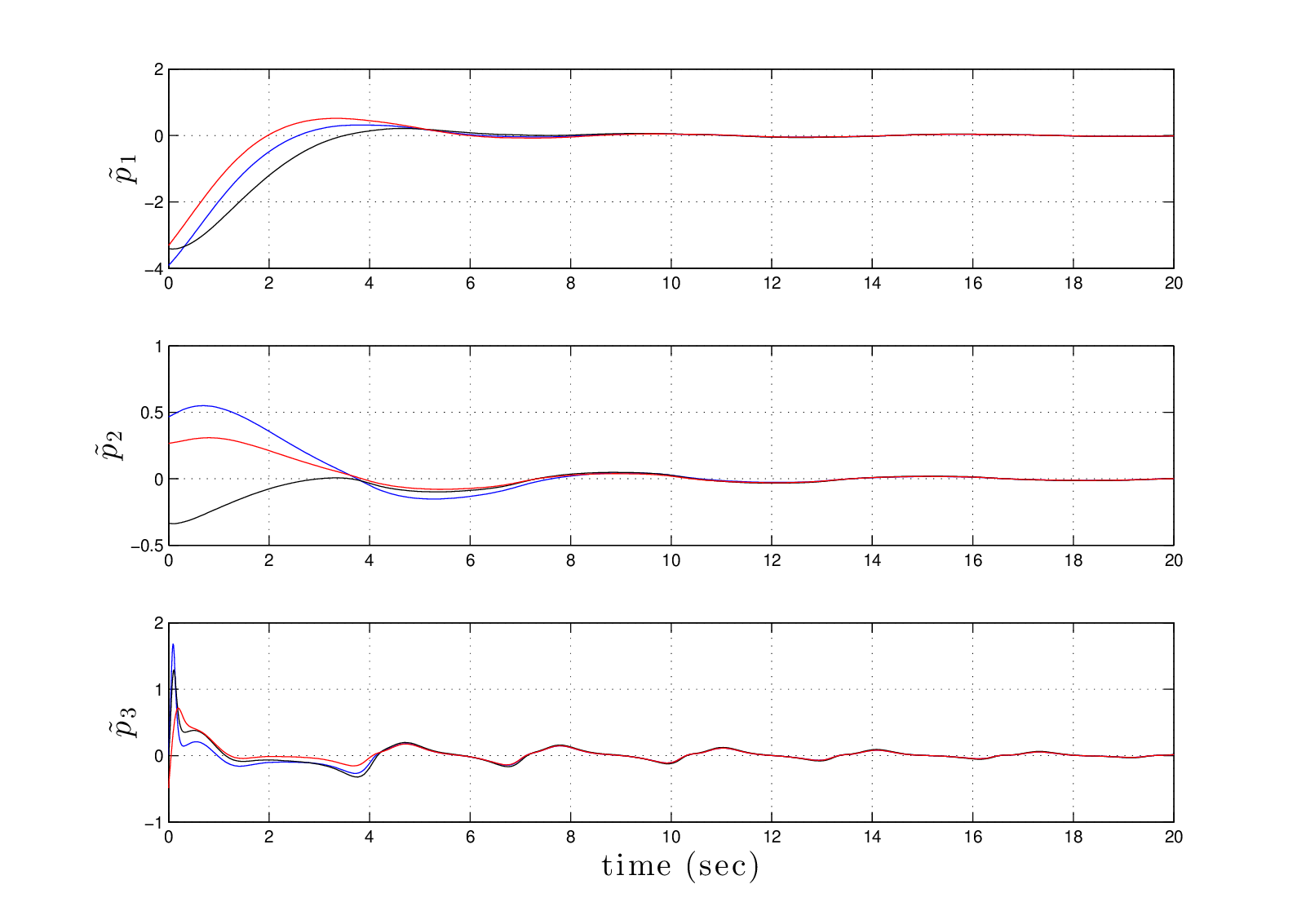}
 \caption{Transient behavior of $\tilde p$ for $\lambda= 0.8$ and different initial conditions of  $d_I,r_{3_I}$ and $P_I$}
 \label{fig2}
\end{figure}

\begin{figure}[h]
 \centering
\includegraphics[width=1\linewidth]{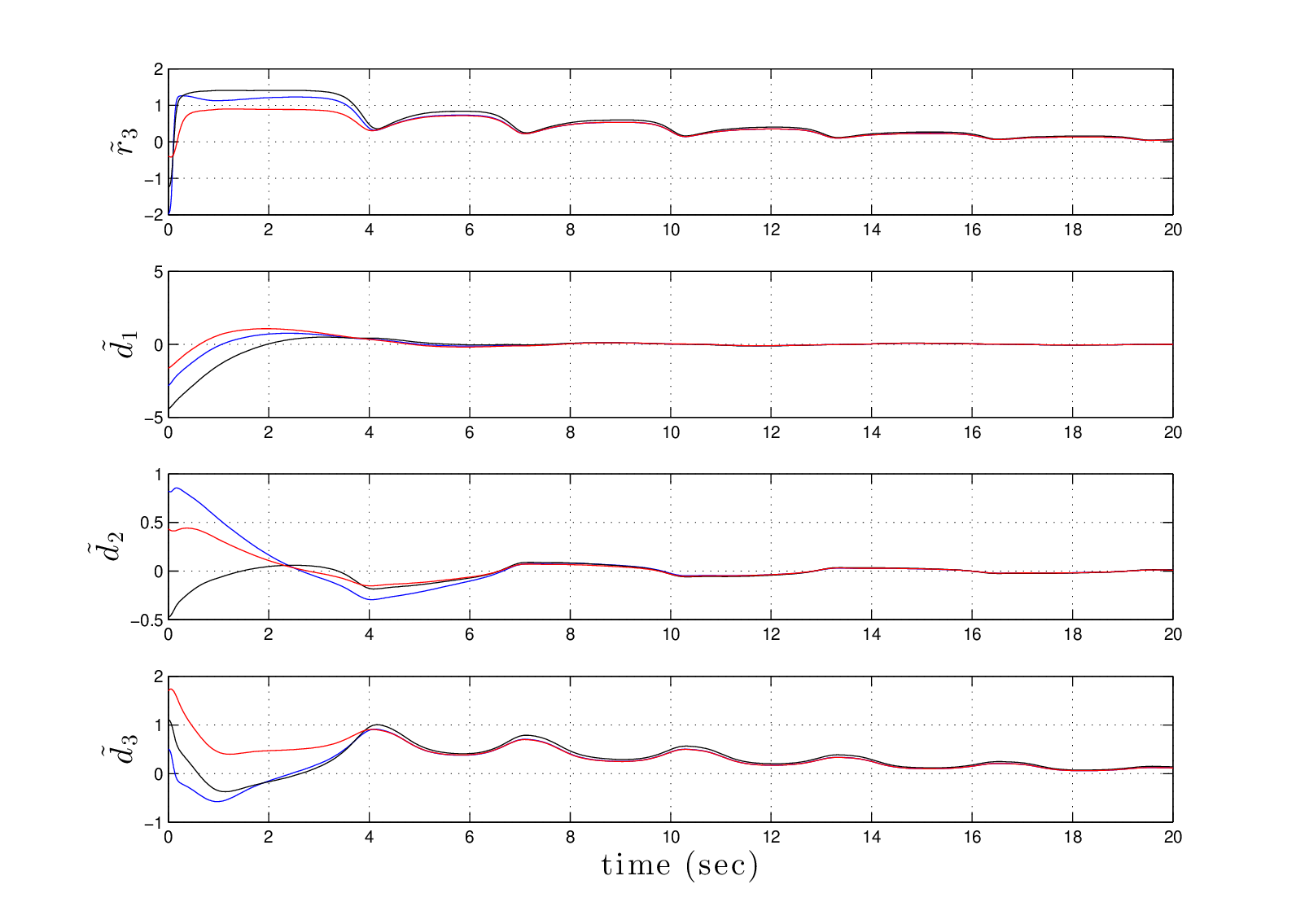}
  \caption{Transient behavior of $\tilde d$ and $\tilde r_3$ for $\lambda= 0.8$ and different initial conditions of  $d_I,r_{3_I}$ and $P_I$}
   \label{fig3}
\end{figure}

The transient behavior of the error signals $\tilde p, \tilde r_3$ and  $\tilde d$ with the tuning parameter $\lambda= 0.8$ and different initial conditions of  $d_I,r_{3_I}$ and $P_I$ are shown in  Fig.\ref{fig2} and Fig.\ref{fig3}. As seen from the figures, besides the convergence to zero of the momenta estimate predicted by the theory, we also observe that the estimated parameters converge to their true value---assessing the fact that the signals chosen for the simulation are persistently exciting.

\begin{figure}[htp]
 \centering
\includegraphics[width=1\linewidth]{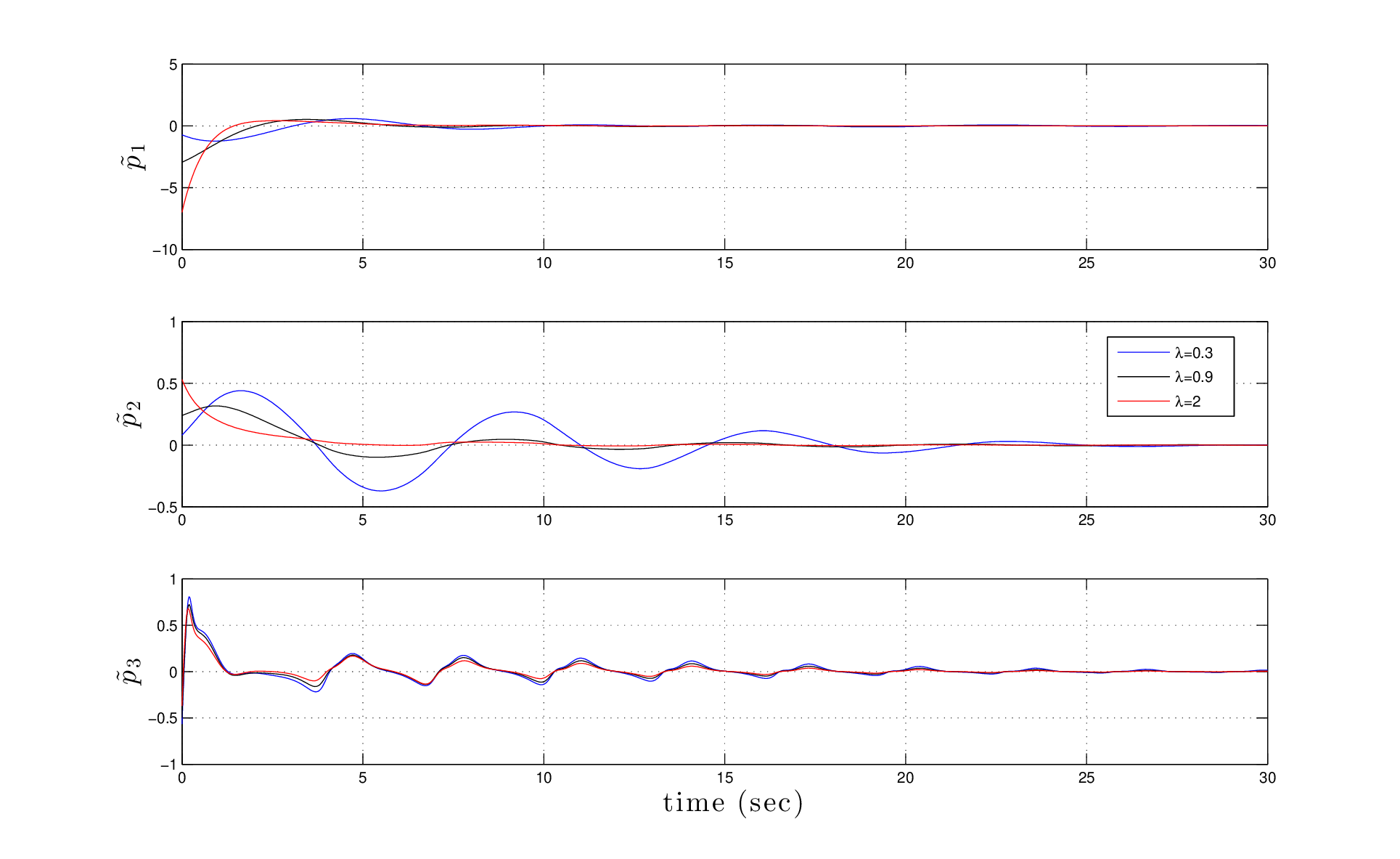}
 \caption{Transient behavior of $\tilde p$ for different values of $\lambda$}
 \label{fig4}
\end{figure}

\begin{figure}[htp]
 \centering
\includegraphics[width=1\linewidth]{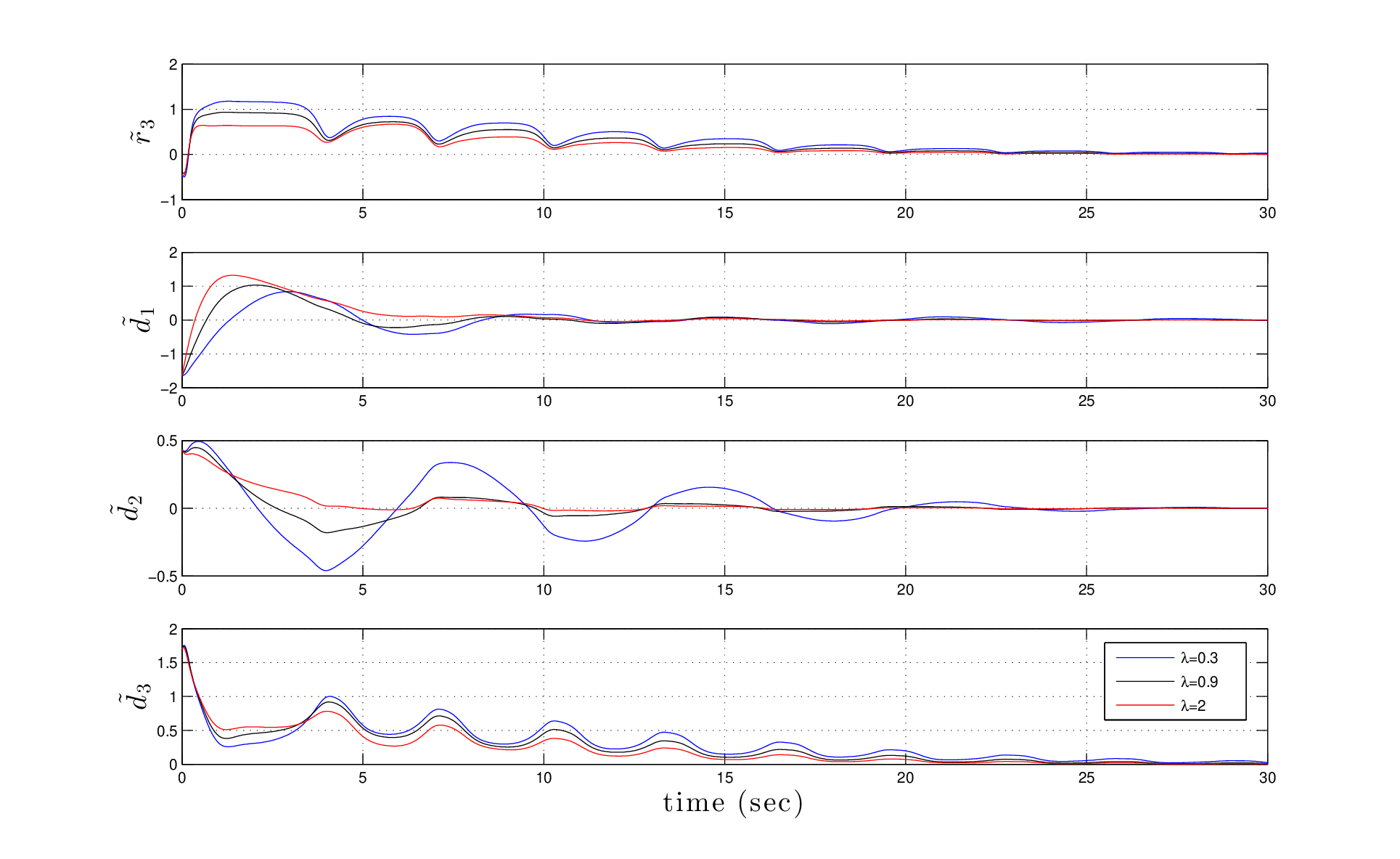}
  \caption{Transient behavior of $\tilde d$ and $\tilde r_3$  for different values of $\lambda$}
   \label{fig5}
\end{figure}

To evaluate the effect of the tuning parameter $\lambda$ on the transient behavior we also  show  in Fig.\ref{fig4} and Fig.\ref{fig5} the transient behavior of the error signals for different values of $\lambda$.
 
\begin{figure}[htp]
 \centering
\includegraphics[width=0.85\linewidth]{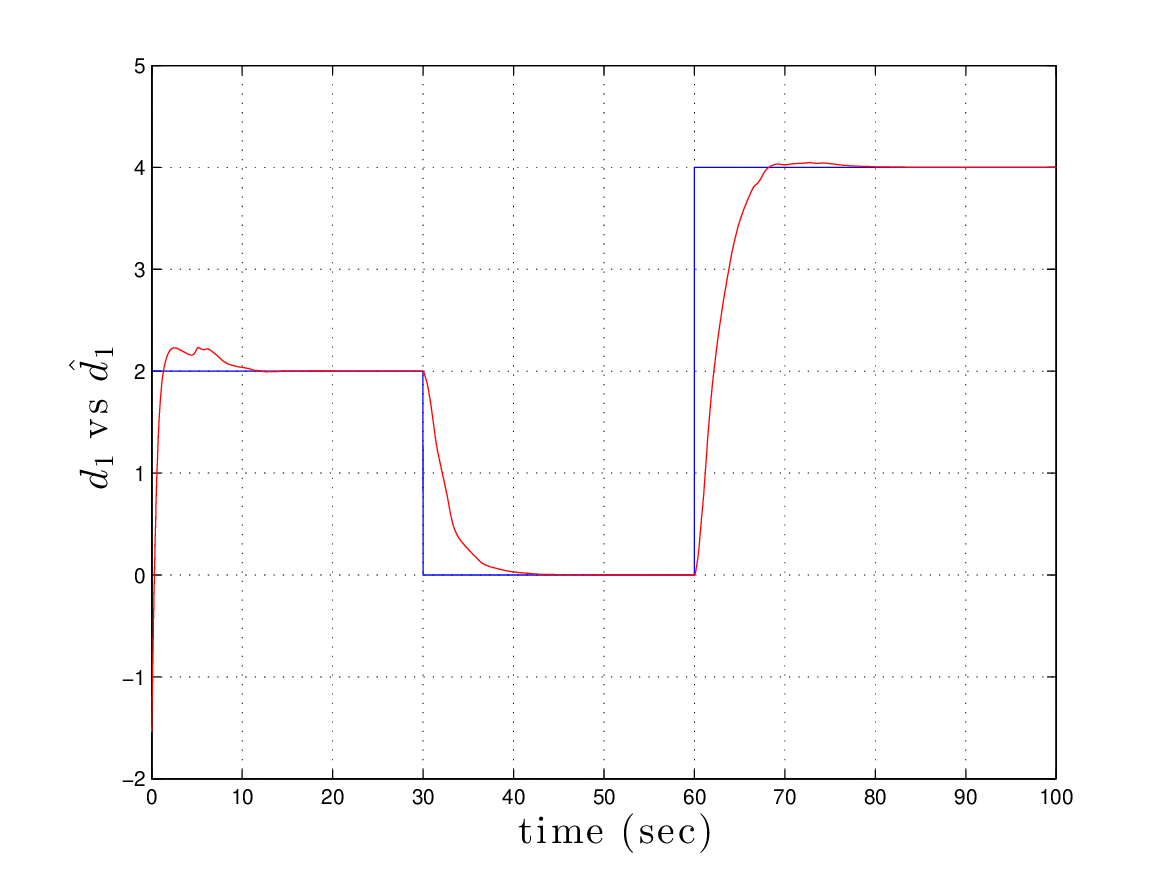}
  \caption{Transient behavior of $d_1$ and $\hat d_1$ with $\lambda=2$}
   \label{fig6}
\end{figure}

\vspace{2mm}
Finally to illustrate the robustness of the adaptive observer,  we carried out  a simulation considering that the input disturbance $d_1$ is subject to step changes.  The  trajectories of $d_1(t)$ and its estimation, depicted in Fig. \ref{fig6}, clearly illustrate the tracking capability of the proposed observer.

\begin{rmrk}\em
It is interesting to note that the standard (lower triangular) Cholesky factorization of  $ M^{-1}(q)$ does not satisfy \eqref{titj}.
\end{rmrk}

\section{Future Research}
\label{sec7}
%
The design of the observer in Proposition \ref{pro1} requires the {\em explicit} solution of the PDE (\ref{qtminone}). This requirement restricts the practical applicability of the approach. Indeed, this PDE
has no {\em free parameters} and its explicit solution may be even impossible. As indicated in \cite{VENetal} this is the case of the classical cart--pendulum example.

Current research is under way to extend the realm of application of the observer in Proposition \ref{pro1}. In particular, it is possible to consider the following generalization of the Lyapunov function
candidate \eqref{V21}
$$
|\tilde p|^2 + \tilde r_u^\top P^{-1} \tilde r_u +|\tilde d|^2,
$$
with $P >0$ a constant matrix. It is straightforward to show that the second PDE in \eqref{PDE_betaa}---that cancels the cross term appearing in the derivative of the new Lyapunov function---becomes
$$
\nabla r_{u_P}  =  -{1 \over \lambda} (\sum_{i=1}^n P Y_i^\top \hat p_i),
$$
where we underscore the presence of the matrix $P$ in front of $Y_i$. There are inertia matrices where the corresponding $Y_i$ are {\em not constant} but there exists positive definite $P$ that will make
$PY_i^\top$ constant---hence relaxing Assumption \ref{ass2}. We are currently investigating whether there exist physical systems for which such property holds.

Another, quite challenging, task is the extension of Proposition \ref{pro1} to systems that do not have ZRS. One possibility is to look into the next class of systems partially linearizable via coordinate changes characterized in \cite{VENetal}.

\begin{ack}                               
This work was supported by the Ministry of Education and Science of Russian Federation (Project 14.Z50.31.0031).
\end{ack}


\end{document}